\documentclass{amsart}

\usepackage{amsmath, amsfonts, amssymb, amscd}

\newtheorem{theo}{Theorem}[section]
\newtheorem{lem}[theo]{Lemma}
\newtheorem{prop}[theo]{Proposition}

\newtheorem{rem}[theo]{Remark}

\newenvironment{pf}{\noindent{\it Proof.
}}{$\blacksquare$\par\medskip}

\newcommand{\C}{{\mathbb C}}
\newcommand{\R}{{\mathbb R}}

\newcommand{\CC}{{\mathcal C}}
\newcommand{\g}{{\mathfrak g}}
\renewcommand{\l}{{\mathfrak l}}
\renewcommand{\k}{{\mathfrak k}}
\newcommand{\m}{{\mathfrak m}}
\newcommand{\z}{{\mathfrak z}}

\newcommand{\h}{{\mathfrak h}}
\newcommand{\p}{{\mathfrak p}}

\renewcommand{\t}{{\mathfrak t}}

\newcommand{\B}{{\mathcal B}}

\newcommand{\U}{{\operatorname {U}}}

\newcommand{\ch}[1]{{#1}^{\vee}}
\renewcommand{\varpi}{\psi}
\newcommand{\tom}{{\tilde\omega}}

\newcommand{\su}{{\mathfrak{su}}}
\newcommand{\so}{{\mathfrak {so}}}
\newcommand{\SU}{{\operatorname{SU}}}
\newcommand{\SO}{{\operatorname{SO}}}

\newcommand{\SL}{{\operatorname{SL}}}
\newcommand{\CP}{{\C\!\operatorname{P}}}
\newcommand{\co}{m}

\newcommand{\Hom}{{\operatorname{Hom}}}

\renewcommand{\span}{\operatorname{span}}

\newcommand{\Ad}{\operatorname{Ad}}
\newcommand{\ad}{\operatorname{ad}}
\newcommand{\reg}{\operatorname{reg}}

\renewcommand{\div}{\operatorname{div}}

\renewcommand{\=}{\overset{\text{def}}{=}}

\def\sideremark#1{\ifvmode\leavevmode\fi\vadjust{
\vbox to0pt{\hbox to 0pt{\hskip\hsize\hskip1em
\vbox{\hsize3cm\tiny\raggedright\pretolerance10000
\noindent #1\hfill}\hss}\vbox to8pt{\vfil}\vss}}}

\title[K\"ahler-Ricci solitons on homogeneous toric bundles (I)]
{K\"ahler-Ricci solitons \\ on homogeneous toric bundles (I)}
\author{Fabio Podest\`a
and Andrea Spiro}


\begin{document}

\begin{abstract} This is the first of a sequence of two papers. Here, a simple algebraic
  characterization of the Fano manifolds in the class of
homogeneous toric bundles over
a flag manifold $G^C/P$ is provided in terms of symplectic
data. The result of this paper is  used in the second paper, where it is proved that an
homogeneous toric bundle over a flag manifold  admits a
K\"ahler-Ricci solitonic metric if and only if it is Fano.
\end{abstract}

\maketitle


\section{Introduction}
\bigskip
This is the first of a sequence of two papers, where  we focus on a particular class of
homogeneous bundles $M$, namely on bundles having a compact toric K\"ahler manifold $F$ as
 fiber and a
generalized flag manifold $G^\C/P=G/K$ as basis, where $G$ is  a
compact
semisimple Lie group, $G^\C$ its complexification and $P$
a suitable parabolic subgroup of $G^\C$. More precisely,
we consider
a surjective homomorphism $\tau:P\to (T^\co)^\C$, where
$T^\co$
is an $\co$-dimensional torus acting effectively on
the toric K\"ahler manifold $F$, $\dim_\C F = m$,
by holomorphic isometries; We then define $M$ to be the
compact complex manifold
$M = G^\C\times_{P,\tau}F$. Any manifold of this kind is
a toric bundle over $G^\C/P$ (see \cite{Ma}) with
holomorphic projection $\pi: M \to G^\C/P$ and it is
  almost $G^\C$-homogeneous (see \cite{HS}) with
  $G$-cohomogeneity equal to $m$. We
  will call any such manifold a {\it homogeneous toric
bundle}. \par
  \smallskip
  The homogeneous toric bundles appear to be direct generalizations of the
$\CP^1$-bundles over flag manifolds studied in
\cite{Sa, KS, KS1, Ko,DW,PS}. These $\CP^1$-bundles are known to
be K\"ahler-Einstein if and only if they are
Fano and their
Futaki functional vanishes identically (\cite{KS,PS});
Moreover
they always admit a K\"ahler metric which is a K\"ahler-Ricci soliton,
provided the first Chern class is positive
(see \cite{Ko,TZ,TZ1}). We also mention that for toric
bundles a uniqueness result for extremal metrics in a
given K\"ahler class is proved in \cite{Gu}. \par
\medskip
Aiming at investigating the existence of K\"ahler-Einstein metrics
and, more generally, of K\"ahler-Ricci solitons in the class of
almost homogeneous toric bundles, in this paper we are interested
in finding simple conditions on $G^\C/P$, $F$ and the homomorphism
$\tau$ which guarantee that the above defined manifold $M$ has
positive first Chern class.   In the next part (\cite{PS2}), we will use such characterization of Fano
homogenous
toric bundles over flag manifolds to show that
any such bundle admits a Ricci-K\"ahler solitonic metric.  As a consequence, in that paper
we obtain
 the following generalization of Koiso and Sakane's theorem:
{\it A homogeneous toric bundle over a flag manifold
is K\"ahler-Einstein if and only if it  is Fano and its Futaki functional
vanishes identically.\/}\par
\medskip
In order to state our characterization  of homogeneous toric bundles with positive first Chern class, we first need to fix some
notations. \par
Let $J$ be the $G^\C$-invariant complex structure on
the flag manifold $G^\C/P=G/K$ and let $\CC$ be the corresponding
positive Weyl chamber in the Lie algebra $\z(\k)$ of the center
$Z(K)$ of $K$ (see e.g. \S 2, for the definition). We will also
use the symbol $\ch {\CC}$ to denote the chamber in $\z(\k)^*$,
which is the image of $\CC$ by means of the dualizing map $X
\mapsto -\B(X,\cdot)$, where $\B$ is the Cartan Killing form on
the Lie algebra $\g$ of $G$.\par
  It is well known that $G^\C/P$ admits a unique
$G$-invariant K\"ahler-Einstein metric $g$
with Einstein constant $c=1$ (see e.g. \cite{Be}): We set
$\mu:G^\C/P\to \g^*$ to be the moment map relative to the
K\"ahler form $\omega =\frac{1}{2 \pi}ÿg(J\cdot,\cdot)$.\par
If $F$ is supposed to be Fano, the Calabi-Yau theorem
implies that for any $T^\co$-invariant K\"ahler form
$\rho \in c_1(F)$, there exists a unique K\"ahler form
$\omega_\rho$ in $c_1(F)$ such that $\rho$ is the Ricci
form
of $\omega_\rho$. In particular, also $\omega_\rho$ is
$T^\co$-invariant. Moreover,
since $b_1(F)=0$ there exists a moment map $\mu_\rho: F
\to \t^*$
relative to $\rho$, uniquely determined up to a constant.
We will say that $\mu_\rho$ is
  {\it metrically normalized\/} if $\int_F \mu_\rho
\omega_\rho^\co = 0$. In \S 4, we will show that the
  convex polytope $\Delta_F = \mu_\rho(F)$, which is the
image of a metrically normalized moment map $\mu_\rho$,
  is actually independent of $\rho$ and it can be
explicitly determined just using the
   $T^\co$ action. \par
Finally, for any given homomorphism
$\tau:P\to (T^\co)^\C$, we set $\mu_\rho^{\tau} \=
(\tau|_{\z(\k)})^*\circ \mu_\rho$. Notice that the map
$\mu^{\tau}_\rho$
is a moment map for the action of
$Z(K)^o$ on $F$ induced by $\tau$, and its image is the
convex polytope $\Delta_{\tau, F} =
  (\tau|_{\z(\k)})^*(\Delta_F) \subset \z(\k)^*$.
Our main result can be now stated as follows.\par
\medskip
\begin{theo} \label{maintheorem} Let $F$ be a toric K\"ahler manifold of
dimension $\co$. Then, for any homomorphism $\tau:P\to
(T^\co)^\C$,
the manifold $M = G^\C\times_{P,\tau}F$ has positive first
Chern class if and only if $F$ is Fano and
$$\mu(eP) + \Delta_{\tau,F}\subset \ch{\CC}\ .\eqno(1.1)$$
\end{theo}
\medskip
Later we will also show that there is a simple algorithm
to determine $\Delta_F$
and that (1.1) can be reformulated in a finite number
of algebraic conditions, which are suited for computations
(see \S 6). \par
\medskip
We mention here that the proof of our main
result originates from an idea for
computing the first Chern class, which goes back to \cite{Ya} and
\cite{Fu}.\par
\medskip
\subsection{Notations} For any Lie group $G$, we will
denote its Lie algebra by the corresponding gothic letter $\g$;
Given a Lie homomorphism $\tau:G\to G'$, we will always use the
same letter to represent the induced Lie algebra homomorphism
$\tau:\g\to\g'$. If $G$ acts on a manifold $M$, for any $X\in
\g$, we will use the symbol $\hat X$ to indicate
  the induced vector field on $M$; We recall here that
$\widehat{[X,Y]} = - [\hat X,\hat Y]$ for every
  $X,Y\in \g$. We denote by $M_{\operatorname{reg}}$ the set of
  $G$-principal points in $M$.
\par
  The Cartan Killing form of a semisimple Lie algebra $\g$
will be always denoted by $\B$ and, for any
  $X\in \g$, we set $\ch X= - \B(X, \cdot) \in \g^*$;
Given a root system $R$ w.r.t to a
fixed maximal torus, we will denote by $E_\alpha\in\g^\C$
the root vector corresponding to the root $\alpha$ in
the Chevalley normalization and by
$H_\alpha=[E_\alpha,E_{-\alpha}]$ the $\B$-dual of
$\alpha$.\par
\bigskip
\bigskip
\section{Preliminaries}
\bigskip
Throughout the following we will denote by
$G$ a connected compact, semisimple Lie group and by
$V=G/K$ a
generalized flag manifold associated to $G$. If we fix a
$G$-invariant complex structure $J_V$ on $V$, then the
complexified
group $G^\C$ acts holomorphically on $V$, which can be
then
represented as $G^\C/P$ for some suitable parabolic
subgroup $P$.\par
We recall that the Lie algebra $\g$ of $G$ admits an
$\Ad(K)$-invariant
  decomposition $\g=\k \oplus \m$ and that,
for any fixed CSA $\h \subset \k^\C$ of $\g^\C$, the
corresponding root system $R$ admits a corresponding decomposition
$R = R_o + R_\m$, so that $E_\alpha \in \k^\C$ if $\alpha \in R_o$
and $E_\alpha \in\m^\C$ if $\alpha\in R_\m$;
Furthermore, $J_V$ induces a splitting $R_\m = R_\m^+\cup R_\m^-$
into two disjoint subset of positive and negative roots, so that
the $J_V$-holomorphic and $J_V$-antiholomorphic subspaces of
$\m^\C$ are given by
$$\m^{(1,0)} = \sum_{\alpha\in R_\m^+}\C E_\alpha,\quad
\m^{(0,1)} = \sum_{\alpha\in R_\m^-}\C
E_\alpha\ .\eqno(2.1)$$
The Lie algebra $\p$ of the parabolic subgroup $P$ is $\p
= \k^\C + \m^{(0,1)}$. It
is also well-known that $G^\C$ is an algebraic group and
$P$ an
algebraic subgroup. \par
Finally, we recall that for any $G$-invariant K\"ahler
form $\omega$ of $V$ there
exists a uniquely associated
  element $Z_\omega \in \z(\k)$ so that
$\left.\omega(\hat X, \hat Y)\right|_{eK} = \B(Z_\omega,
[X,Y])$
for any $X, Y\in \g$. Notice that, by (2.1) and the
positivity of $\omega$, the element
$Z_\omega$ has to belong to the positive Weil chamber
$$\CC = \{\ W \in \z(\k)\ :\ i\alpha(W) > 0\ ,\ \text{for
any}\ \alpha \in R^+_\m\ \}\ .$$
  Moreover,
a straightforward check shows that the moment map
$\mu_\omega: V \to \g^*$ relative to $\omega$ is given by
$\mu_\omega(g K) = (\Ad_gZ_\omega)^\vee$ for any $gK \in
V$.
We recall also that the
$G$-invariant K\"ahler-Einstein form $\omega_V$ on $V$
with Einstein constant $c =1$, is determined by the
associated element
(see e.g. \cite{Be, BFR} - be aware that in this paper,
we adopt the definition of Ricci form $\rho$ used e.g. in \cite{Fu1}, which differs from
the one in  \cite{Be} and
\cite{BFR}   by the factor $\frac{1}{2 \pi}$)
$$Z_V = - \frac{1}{2 \pi} \sum_{\alpha \in R^+_\m} i
H_{\alpha}\ .\eqno(2.2)$$
We will now focus on those flag manifolds
$G^\C/P =G/K$ for which there exists a surjective
homomorphism $\tau:P\to (T^\co)^\C$. Using the structure
of parabolic subgroups and the fact that $(T^\co)^\C$ is
abelian, we
see that $\tau$ is completely determined by its
restriction to $K$; Moreover $\tau|_K$ takes value in
$T^\co$ and hence it is fully determined by
its restriction to
the connected component $Z^o(K)$ of the center of the
isotropy $K$. We can
therefore
consider the algebraic manifold
$$M \= G^\C \times_{P,\tau} F = G\times_{K,\tau} F\ ,$$
where $P$ (or $K$) acts on $F$ by means of
$\tau$.\par
\medskip
  The manifold $M$ is a fiber
bundle
over the flag manifold $G/K$ with holomorphic projection
$\pi$; Moreover $G^\C$ acts
almost homogeneously, i.e. with an open and dense orbit in
$M$, while $G$ acts
by cohomogeneity $\co$ with principal isotropy type $(L)$,
where $L = \ker \tau\cap G\subset K$. \par
\smallskip
We prove now the following Lemma, which will be useful and
often
tacitly used in the sequel.
\begin{lem} \label{complexstructure} If $J$ denotes the
complex structure of $M$, then
$$J\hat \m_p = \hat \m_p$$
for every $p\in \pi^{-1}([eK])$.
\end{lem}
\begin{pf} We know that $\m^{(0,1)}\subset \p$ and
$\tau^\C |_{\m^{(0,1)}}=0$, so that
$\hat \m^{(0,1)}|_p = 0$. The vectors
$E_\alpha -E_{-\alpha}$ and $i(E_\alpha + E_{-\alpha})$,
$\alpha \in R^+_\m$, span $\m$ over the
reals
and
$$J(\widehat{(E_\alpha -E_{-\alpha})}|_p) =
J(\hat{E_\alpha}|_p) =
\widehat{iE_\alpha}|_p =
\widehat{i(E_\alpha + E_{-\alpha})}|_p\in \hat \m_p\ .$$
Similarly, $J(\widehat{i(E_\alpha + E_{-\alpha})})\in\hat
\m_p$.
\end{pf}
\smallskip
In the sequel $(Z_1,\ldots,Z_\co)$ will denote a fixed
$\B$-orthonormal basis of $(\ker \tau)^\perp \cap
\z(\k)$.\par
\smallskip
\bigskip
\bigskip
\section{The algebraic representatives}
\bigskip
We recall that, if $\varpi$ is a $G$-invariant 2-form on $M$,
for any $p\in M$
there exists a unique
$\operatorname{ad}_{\g_p}$-invariant
element $F_{\varpi, p}\in \Hom(\g,\g)$ such that
$\B(F_{\varpi, p}(X), Y) = \varpi_{p}(\hat X, \hat Y)$ for
any $X, Y \in \g$.
Moreover, if $\varpi$ is closed, it turns out that
$F_{\varpi, p}$ is a derivation of $\g$ and hence of the
form
$\ad_{Z_\varpi}$ for some element $Z_\varpi$ belonging to
the centralizer in $\g$ of the
isotropy subalgebra $\g_p$ (see e.g.
\cite{PS,PS1,Sp}).\par
So, in the following, for any $G$-invariant, closed
2-form $\varpi$, we will denote by $Z_\varpi$
the $G$-equivariant map $Z_\varpi: M \to \g$ defined by
$$\varpi_p(\hat X, \hat Y) = \B([Z_\varpi|_p, X], Y) =
\B(Z_\varpi|_p, [X, Y])\qquad \text{for any}\ X, Y \in
\g\eqno(3.1)$$
and it will be called {\it algebraic representative of\/}
$\varpi$.\par
Notice that, in case $\varpi$ is non-degenerate,
the moment map $\mu_\varpi: M \to \g^*$ relative to
$\varpi$
coincides with the $(-\B)$-dual map of the algebraic
representative, i.e. $\mu_\varpi = \ch Z_\varpi$.
In fact, by the closure and $G$-invariance of $\varpi$, we
have that for any
vector field $W$ on $M$ and any $X, Y\in \g$
$$0 = d \varpi(W, \hat X, \hat Y) = W(\varpi(\hat X, \hat
Y)) +
\varpi(W, [\hat Y, \hat X]) = $$
$$ = W(\varpi(\hat X, \hat Y)) +
\varpi(W, \widehat{[X, Y]}) = \B(W(Z_\varpi), [X,Y]) +
\varpi(W, \widehat{[X, Y]}) \ .
\eqno(3.2)$$
Since $\g = [\g, \g]$, it follows that $d(\ch
Z_\varpi)(W)(X) = \varpi(W,\hat X)$,
which implies the claim.\par
\medskip
By $G$-equivariance, notice that any algebraic
representative $Z_\varpi$ is uniquely determined
by its restriction on the fiber $F = \pi^{-1}(eK) \subset
M$. Such restriction satisfies the following.
\par
\begin{lem} \label{recovering}
Let $\varpi$ be a $G$-invariant, $J$-invariant, closed
2-form on $M$.
\begin{itemize}
\item[(a)] If the restriction $\psi|_F$ satisfies
$\psi(\hat Z_j,\hat \m)=0$ for $1\leq j\leq\co$, then
$Z_\psi|_F$ is of the form $Z_\varpi|_F = \sum_{i = 1}^\co
f_i Z_i + I_\varpi$, where $I_\varpi \in \l = Lie(\ker \tau\cap G)$ and
$f_i: F
\to \R$ are smooth functions. Moreover $\varpi_p(J \hat
Z_j, \hat Z_i) = \left. J \hat Z_j(f_i)\right|_p$ for any
$p \in F$ and $1 \leq i,j \leq \co$ and $I_\varpi$ is
constant;
in particular, $\varpi$ can be completely recovered by
its algebraic representative $Z_\varpi$ (and hence
by its associated moment map, if $\varpi$ is
non-degenerate);
\item[(b)] $\varpi$ is cohomologous to $0$ if and only if
$Z_\psi|_F = - \sum_i J \hat Z_i(\phi)Z_i$
for some $K$-invariant smooth function $\phi: F \to \R$
and, if this occurs, then $\psi = d d^c \phi$.
\end{itemize}
\end{lem}
\begin{pf} (a)\ Since $\k = \l
+\span\{Z_1,\ldots,Z_\co\}$, $[\k,\m] = \m$ and $\hat
\l|_F = 0$, we have that on $F$
$$0 = \psi(\hat \k, \hat \m) = \B(Z_\psi,[\k,\m]) =
\B(Z_\psi,\m)\ ;$$
hence $Z_\psi|_F$ takes values in $\k$ and has the
claimed form. From (3.2) we have that for any $X, Y\in
\g$ and any $p\in F$,
$$\varpi_p(J \hat Z_i, \widehat{[X,Y]}) = -
\B\left(\sum_{j=1}^\co
J \hat Z_i(f_j)|_p Z_j + J \hat Z_i( I_\varpi)|_p, [X,
Y]\right)\ .\eqno(3.3)$$
On the other hand,
$$\varpi_p(J \hat Z_i, \widehat{[X,Y]})
= - \sum_{j = 1}^\co \B([X,Y], Z_j) \varpi_p(J \hat Z_i,
\hat Z_j) -
\varpi_p(\hat Z_i, J \widehat{[X,Y]_{\m}})\ ,\eqno(3.4)$$
where we denote by $(\cdot)_{\m}$ the $\B$-orthogonal
projection onto $\m$. Since $\hat \m_p$ is $J$-invariant
and $\varpi_p$-orthogonal to
$\span\{\hat Z_i|_p\}$, the second term of (3.4) vanishes and,
from (3.3), we obtain
$$\sum_{j = 1}^\co \B\left(\varpi_p(J \hat Z_i, \hat Z_j)
Z_j,[X,Y]\right) =
\B\left(\sum_{j=1}^\co
J \hat Z_i(f_j)|_p Z_j + J \hat Z_i( I_\varpi)|_p, [X,
Y]\right)\ .$$
Since $[\g,\g]=\g$, we have that
$$\varpi(J \hat Z_i, \hat Z_j) \equiv J \hat
Z_i(f_j)\ ,\qquad J \hat Z_i( I_\varpi) = 0\ .\eqno(3.5)$$
This together
with the fact that
   $\hat Z_i(I_\varpi) = 0$, which is due to
$G$-equivariance, implies that $I_\varpi$ is
constant on $F$ and the first claim follows. Furthermore, formula (3.3)
implies that if the map
$Z_\varpi|_F = \sum_{i=1}^\co f_i Z_i + I_\varpi: F \to
\z(\k)$ is known,
  it is possible to evaluate $\varpi_p(\hat X, \hat Y)$
for any $X, Y\in \g^\C$ and $p\in F$; then,
from almost homogeneity also the last claim of (a) follows. \par
\noindent (b)\ Notice that in case $\varpi$ is
cohomologous to $0$, then
it is of the form
$\varpi = dd^c \phi$ for some $G$-invariant real valued
function $\phi$. Then, for any $X, Y\in \m$, on $F$ we have that
$$d d^c \phi(\hat X, \hat Y) = - \hat X(J \hat Y(\phi)) +
\hat Y(J \hat X(\phi)) + J [\hat X, \hat Y](\phi) = $$
$$ = J \widehat{[X,Y]}(\phi) = - \sum_{i=1}^\co \B(Z_i,
[X,Y]) J\hat Z_i(\phi)\ .$$
It follows immediately that $Z_\varpi|_F = - \sum_{i=1}^\co
J\hat Z_i(\phi) Z_i$. Conversely,
if $Z_\varpi|_F= -
\sum_{i=1}^\co J\hat Z_i(\phi) Z_i$
  for some $K$-invariant function $\phi\in C^\infty(F)$,
then by (a) and the above remarks, $\varpi=d d^c \phi$,
where we consider $\phi$ as $G$-invariantly extended to
$M$.
\end{pf}
\bigskip
  We want now to determine the algebraic representative of
the Ricci form
$\rho$ of a given $G$-invariant K\"ahler form $\omega$. In
what follows, we denote by $(F_\alpha,
G_\alpha)_{\alpha\in R^+_\m}$
the basis for $\m$ given by the vectors
$$F_\alpha = \frac{1}{\sqrt{2}} \left(E_\alpha -
E_{-\alpha} \right)\ ,\qquad
G_\alpha = \frac{i}{\sqrt{2}} \left(E_\alpha + E_{-\alpha}
\right)$$
with $\alpha \in R^+_\m$. Notice that, by definition of
$R^+_\m$,
$J_V F_{\alpha} = G_\alpha$ and $J_V G_\alpha = -
F_\alpha$
and that the complex structure $J$ of $M$ induces on $\m$
the same complex structure
induced by $J_V
$ (see proof of Lemma \ref{complexstructure}).
We order the roots in $R^+_\m$ so to call them $\alpha_1$,
$\alpha_2$, etc.,
and we denote by $F_i$, $1 \leq i \leq \co + |R^+_\m|$ the
elements of $\k$ defined by
$F_i = Z_i$ if $1\leq i \leq \co$ and $F_i =
F_{\alpha_{i-\co}}$ if $\co + 1 \leq i$.
Notice that, at any point $p\in F \cap M_{\reg}$ the
vector fields $\{\hat F_j,
J \hat F_j\}$ are linearly independent and span the whole
$T_p M$. Finally,
for any given K\"ahler form $\omega$, let us also denote
by
$h: M \to \R$ the function
$$h(q) = \omega^n(\hat F_1, J \hat F_1, \hat F_2,
J \hat F_2, \dots )|_q\ .\eqno(3.6)$$
\medskip
We claim that, at any point $p \in F$, $\hat X(h)|_p = 0$
for any $X\in \g$. In fact, using ${\mathcal L}_{\hat X}\omega = 0$
and ${\mathcal L}_{\hat X}J = 0$, we have
$$\hat X(h)|_p = - \omega^n(\widehat{[X,F_1]}, J \hat
F_1, \hat F_2,
J \hat F_2, \dots ) - \omega^n(\hat F_1, J \widehat{[X,
F_1]}, \hat F_2,
J \hat F_2, \dots ) - $$
$$ - \omega^n(\hat F_1, J \hat F_1, \widehat{[X, F_2]},
J \hat F_2, \dots ) - \dots\ .$$
On the other hand, for any $i$,
$$[X, F_i] \in \l + \span\{ F_j,\ j \neq i\ \} +
\span\{ J_V F_j,\ j\geq \co +1\ \}\ .$$
and this implies that
$$\omega^n(\hat F_1, \dots, J F_{j-1}, \widehat{[X,F_j]},
J F_j, \dots) =
\omega^n(\hat F_1, \dots, J F_{j-1}, F_j, J \widehat{[X,F_j]},
 \dots) = 0$$
and hence the claim. We may now prove the following.
\begin{prop}\label{ricci} The restriction to $F_{\reg} = F\cap M_{\reg}$
of the
algebraic representative $Z_\rho$ of the Ricci form $\rho$
of a K\"ahler form $\omega$ is
$$Z_\rho = \sum_{i = 1}^\co \frac{J \hat Z_i(\log |h|)}{4 \pi}
Z_i +
Z_V\ ,\eqno(3.7)$$
where $h$ is the function (3.6) and $Z_V \in \z(\k)$ is
the element defined in (2.2).
Furthermore, if $Z_\omega = \sum_i f_i Z_i + I_\omega$ is
the restriction to $F$ of
  the algebraic representative of $\omega$, then the
function $h$ is
  $$h = K\cdot \det\left(\begin{matrix}
f_{i,j}\end{matrix}\right)\cdot \prod_{\alpha \in R^+_\m}
(\sum_{i=1}^\co a^i_{\alpha} f_i + b_\alpha)\ ,\eqno(3.8)$$
where $f_{i,j} \= J \hat Z_j(f_i)$,
$a^i_{\alpha} \= \alpha(i Z_i)$,
  $b_\alpha \= \alpha(i I_\omega)$ and $K$ is a real
constant.
\end{prop}
\begin{pf}
We first show that $\rho|_F(\hat Z_j,\hat\m)=0$. By
Koszul formula (see e.g. \cite{Be}, p. 89 - be aware
of the difference recalled in \S 2 between the  definitions of $\rho$ adopted in  \cite{Be} and in this paper)
$$\rho_p(\hat X, \hat Y) =
-\frac{1}{4 \pi} \frac{{\mathcal L}_{J\widehat{[X,Y]}}
\omega^n (\hat F_1, J \hat F_1, \hat F_2,
J \hat F_2, \dots )}{
\omega^n (\hat F_1, J \hat F_1, \hat F_2,
J \hat F_2, \dots )}\eqno(3.9)$$
for every $p\in F_{\operatorname{reg}}$ and $X,Y\in \g$. On the other
hand, we
claim
that for any $W\in\m$,
$${\mathcal L}_{J \hat{W}} \omega^n|_F = 0\ .\eqno(3.10)$$
In fact,
$$\B([W, F_{\alpha_k}], G_{\alpha_k}) =
\B(W, [F_{\alpha_k}, G_{\alpha_k}]) = \B(W, i
H_{\alpha_k}) = 0\ .$$
So, for any $\hat F_j$ with $j \geq \co +1$,
$\widehat{[W, F_j]}$
has trivial component along $J \hat F_j$ and
$J\left(\widehat{[W, F_j]}\right)$
has trivial component along
$\hat F_j$. This implies that
$$\omega^n(\hat F_1, \dots, J [\hat W, \hat F_j], \dots )
=
\omega^n(\hat F_1, \dots, J [\hat W, J \hat F_j], \dots )
=
  0\ .$$
Moreover, when $1 \leq i \leq m$,
we have that $[W, F_i] = [W, Z_i]\in \m$ and hence also
in this case
$$\omega^n( J [\hat W, \hat F_1], J \hat F_1, \hat F_2,
\dots) =
\omega^n( \hat F_1, J [\hat{W}, J \hat F_1], \hat F_2,
  \dots) = \dots = 0\ .$$
These facts imply that
$${\mathcal L}_{J\hat{W}} \omega^n (\hat F_1, J \hat F_1,
\hat F_2,
J \hat F_2, \dots ) = J\hat{W}(h) \ .$$
Using $J \hat \m_p = \hat \m_p$ and the fact that $\hat
X(h)_p = 0$
for any $X\in \g$, we get (3.10). \par
Now, the fact that $\rho|_F(\hat Z_j,\hat\m)=0$ follows
immediately from (3.9), (3.10) and
from $[\k,\m]\subseteq \m$.\par
\medskip
By Lemma \ref{recovering}(a) and (3.9), in order to
determine $Z_\rho$,
it suffices to compute for $p\in F_{\operatorname{reg}}$
   $$\rho_p(\hat X, J \hat X) = \rho_p(\hat X, \widehat{J
X}) =\left.
-\frac{1}{4 \pi} \frac{{\mathcal L}_{J\widehat{[X,JX]}}
\omega^n (\hat F_1, J \hat F_1, \hat F_2,
J \hat F_2, \dots )}{
\omega^n (\hat F_1, J \hat F_1, \hat F_2,
J \hat F_2, \dots )}\right|_p\ .\eqno(3.11)$$
Now, given $X\in \m$, we put $Y = [X,JX]$ and we write
$Y = Y_\m + \sum_{i = 1}^\co Y_i Z_i + Y_\l$,
where subscripts indicate $\B$-orthogonal projections onto
the corresponding subspaces.
Using (3.10) and $\hat \l|_F = 0$, (3.11) reduces to
$$\rho_p(\hat X, J \hat X) = - \frac{1}{4 \pi h} \sum_{i =
1}^\co Y_i \left(J \hat Z_i(h)\right)
  + \frac{1}{4 \pi h} \omega^n(J [\widehat{Y_{\k}}, \hat
F_1],
J \hat F_1, \hat F_2, ....) + $$
$$ + \frac{1}{4 \pi h}
\omega^n(\hat F_1, J[ \widehat{Y_{\k}}, \hat F_1],
\hat F_2, ....) + \dots\ .\eqno(3.12)$$
Now, recall that $[\k, F_i] = 0$ for all $1 \leq i \leq
\co$ and that,
  for any $1 \leq j$ and any
$ H\in \h\cap\g$, where $\h$ is the fixed CSA,
$$J \widehat{[H, F_{j+\co}] } =J \widehat{[H, F_{\alpha_{j}}]}|_p =
  - i\alpha_{j}(H) J \widehat{G_{\alpha_{j}}}|_p = i
\alpha_{j}(H) \hat F_{\alpha_j}|_p\ . $$
At the same time, for any $W \in
\span_{\R}\{ E_\beta,\ \beta \in R\}\cap \g$,
the bracket $[W, F_{\alpha_j}]$ is always orthogonal to
$\span_\R\{F_{\alpha_{j}}, G_{\alpha_{j}}\}$. Therefore
$\widehat{[W, F_{\alpha_j}]}|_p$
has trivial component along the vector $\hat
G_{\alpha_j}|_p = J \hat F_{\alpha_j}|_p$ and
$J \widehat{[W, F_{\alpha_j}]}|_p$
has trivial component along the vector $\hat
F_{\alpha_j}$. It follows that
$$ \frac{1}{4 \pi h} \omega^n(J [\widehat{Y_{\k}}, \hat
F_1],
J \hat F_1, \hat F_2, ....) +
\frac{1}{4 \pi h}
\omega^n(\hat F_1, J [\widehat{Y_{\k}}, J \hat F_1],
\hat F_2, ....) + \dots = $$
$$ = -\frac{1}{4 \pi h} \omega^n(\hat F_1, \ldots,
J \hat F_\co, \widehat{[Y_{\h}, F_{\co+1}]},
J \hat F_{\co+1},
  \hat F_{\co+2}, \ldots) - $$
$$ - \frac{1}{4 \pi h} \omega^n(\hat F_1, \ldots,
J \hat F_\co, \hat F_{\co+1}, J \widehat{[Y_{\h},
F_{\co+1}]},
  \hat F_{\co+2}, ....) - \ldots = $$
$$ = - \frac{1}{2 \pi}  \sum_{\alpha\in R^+_\m} i \alpha(Y_{\h}) =
\B\left( Z_V, [X, JX]\right)\ .$$
So, (3.12) reduces to
$\rho_p(\hat X, J \hat X) = \B\left( \sum_{i = 1}^\co
\frac{J \hat Z_i(h)}{4 \pi h} Z_i + Z_V, [X, JX]\right)$ and
(3.7) follows.\par
To check (3.8), it suffices to observe that for any $j
\geq 1$,
$$\omega_p(\hat F_{j+\co}, J \hat F_{j+\co}) = \omega_p(\hat F_{\alpha_j},
\hat
G_{\alpha_j}) =
\B(\sum_i f_i(p) Z_i + I_\omega, [F_{\alpha_j}, G_{\alpha_j}]) = $$
$$ = \B(\sum_i f_i(p) Z_i + I_\omega, i
H_{\alpha_{j}}) = \sum_i f_i(p)
a^i_{\alpha_{j}} + b_{\alpha_{j}} $$
and that
$h(p) = \omega^{\co}_p(\hat F_1, J\hat F_1, \dots, \hat
F_\co, J \hat F_\co) \cdot
\prod_{\co +1\leq j} \omega_p(\hat F_j, J \hat F_j)$.
Then the conclusion follows from Lemma \ref{recovering}
(a).\end{pf}
\bigskip
\bigskip
\section{Canonical polytope of a Fano toric manifold}
\bigskip
In the following, let $F$ be Fano and, as considered in
the Introduction, for any $T^\co$-invariant K\"ahler
form $\rho \in c_1(F)$, let us denote by $\omega_\rho$ the
unique K\"ahler form in
$c_1(F)$ which has $\rho$ as Ricci form. In particular,
also $\omega_\rho$ is $T^\co$-invariant.
Since $b_1(F) = 0$, there is a moment map $\mu_\rho: F \to
\t^*$ which is
uniquely determined up to a constant. We say that
$\mu_\rho$ is {\it metrically normalized\/}
if $\int_F \mu_\rho \omega_\rho^\co = 0$.\par
We now want to show that the convex polytope
$\Delta_{F,\rho} \subset \t^*$, which is image of $F$
under
the metrically normalized moment map $\mu_\rho$, is
actually independent of the choice of $\rho$
and it is canonically associated with $F$. For any
K\"ahler form $\omega$ on $F$, we may construct the map
$\delta_\omega:F\to \t^*$ defined by
$$\delta_\omega|_p (W)\= \frac{1}{4 \pi} \div_p(J\hat W),\qquad
W\in \t,\quad p\in F\ ,\eqno(4.1)$$
where the divergence is defined by ${\mathcal L}_X
\omega^n =
\div(X) \omega^n$.
Notice that the map (4.1) is well defined even when $F$ is
not Fano. Moreover, by standard facts on divergences, we
have that
for any point $p \in Fix(T^\co)\subset F$ and $Z\in \t$
$$\delta_\omega|_p(Z) = \frac{1}{4 \pi} Tr(J\circ A_Z|_p)\ ,\eqno(4.2)$$
  where the $A_Z|_p$ is the linear map $A_Z|_p: T_p F \to
T_p F$ defined by
$$A_Z|_p(v) = \left. \frac{d}{dt}\left(\Phi^{ \hat
Z}_{t*}(v)\right)\right|_{t = 0}\ ,\qquad v\in T_pF\ ,
\eqno(4.3)$$
where $\Phi^{\hat Z}_{t}$ is the flow generated by $\hat
Z$.
In particular $ \delta_\omega(Fix(T^\co))$ does not depend
on
  $\omega$ and it is uniquely determined
  just by the holomorphic action of $T^\co$ on $F$. In
the following,
  we will call the convex hull $\Delta_F \subset \t^*$
  of the points of $ \delta_\omega(Fix(T^\co))$
    the {\it canonical polytope of $(F, T^\co)$\/}.\par
Let us now go back to a $T^\co$-invariant K\"ahler form
$\rho \in c_1(F)$ and to the
K\"ahler form $\omega_\rho$, which has $\rho$ as Ricci
form.
If we denote by $g_{\alpha,\bar\beta}$ the components of
the K\"ahler metric
$g = \omega_\rho(\cdot, J\cdot)$ in a system of
holomorphic coordinates, then the Ricci
form $\rho$ of $\omega_\rho$ is equal to
$$\rho =
-{\frac{1}{4 \pi}}dd^c\log(\det(g_{\alpha,\bar\beta}))$$
and at any $T^\co$-regular point we have
$\det(g_{\alpha,\bar\beta}) = f\cdot h_\rho$,
where $f$ is the squared norm of a holomorphic function and
$h_\rho = \omega_\rho^\co(\hat Z_1,J\hat Z_1,\ldots,\hat
Z_\co,J\hat Z_\co)$. Using
the
fact that $\mathcal L_{\hat Z_i}\omega_\rho^\co = 0$,
we see that
$$\rho(\hat Z_i,J\hat Z_k) = -\frac{1}{4 \pi }J\hat
Z_k\left(\frac{J\hat Z_i(h_\rho)}{h_\rho}\right)
= -\frac{1}{4 \pi}J\hat Z_k\left( \div(J\hat Z_i)\right).\eqno(4.4)$$
From (4.4) and Stokes theorem, one
may check that the map (4.1) with $\omega =
\omega_\rho$
  coincides with the metrically normalized moment map
$\mu_\rho$ relative
to $\rho$. From the previous remarks
it follows immediately that $\Delta_{F,\rho} \=
\mu_\rho(F)$ coincides with the
canonical polytope $\Delta_F$ of $F$ and hence it is
independent of the chosen K\"ahler form $\rho \in
c_1(F)$.\par
\begin{rem} {\rm Notice that when $\omega$ is a
$T^\co$-invariant K\"ahler-Einstein form with Ricci form
$\rho=\omega$, the metrically normalized moment map
$\mu_\rho$ satisfies
$$\int_F\mu_\rho \rho^\co = \int_F \mu_\rho \omega^\co = 0\ ,$$
so that the barycenter of $\Delta_F$ is the origin. This
obstruction to the existence of K\"ahler-Einstein metrics
is known to be equivalent to the vanishing of the Futaki
invariant (see \cite{Ma, Fu1}).
To the best of our knowledge, our characterization of
$\Delta_F$ in terms of the $T^\co$-action is
new.}\end{rem}
\bigskip
\bigskip
\section{Proof of the main theorem}
\bigskip
Let us fix a $G$-invariant K\"ahler form $\omega$ on $G/K$
and let $Z_\omega \in \z(\k)$ be
its associated element. Pick also a $T^\co$-invariant
2-form $\rho \in c_1(F)$ and let $\omega_\rho$ be the
unique $T^n$-invariant
K\"ahler form in a fixed K\"ahler class, which has $\rho$
as Ricci form. Denote also by $\mu_{\omega_\rho}$
a fixed moment map relative to $\omega_\rho$.\par
We now fix the fiber $\pi^{-1}(eP)\cong F$ and define a
2-form $\tom$ on
$TM|_F$ as
follows: for $p\in F$, $X,Y\in T_pF$ and $A,B\in \m$
$$\tom_p(X,Y) = \omega_\rho|_p(X,Y)\ ;\ \quad
\tom_p(X,\hat A) = 0; \ \quad \tom_p(\hat A,\hat
B) = -
\mu_{\omega_\rho}(p)(\tau([A,B]_\k))\ ,\eqno(5.1)$$
where for every $U\in \g$ we denote by $U_\k$ the
component along $\k$
w.r.t. the
decomposition $\g = \k \oplus \m$.
\par
We now extend $\tom$ to a global $G$-invariant 2-form,
still denoted by
$\tom$. We can easily check that $\tom$ is $J$-invariant,
using the fact that for any $A,B\in\m$
we have $[A,B]_{\z(\k)} = [J_VA,J_VB]_{\z(\k)}$. We claim
that
$\tom$ is also closed. It is enough to check that
$d\tom_p(\hat A,X,Y) =
d\tom_p(\hat A,\hat B,Y) = 0$ for any $p\in F$, since the
condition $d\tom_p(\hat A,\hat B, \hat C)=0$ for
$A,B,C\in\m$ follows immediately from the
Jacobi identities in $\g$.
Since for any pair of vector fields $V_1$, $V_2$,
$$0 = \mathcal L_{\hat A}\tom(V_1,V_2) =
d\left(\imath_{\hat A} \tom\right)(V_1,V_2) +
  d\tom(\hat A, V_1,V_2)\ ,$$
we are reduced to check that $\left.d\left(\imath_{\hat
A}\tom\right)(X,Y)\right|_p = 0$ and $\left.d\left(\imath_{\hat A}
\tom\right)(\hat B,Y) \right|_p = 0$ for every $X,Y\in T_pF$. Now, if we
extend $X,Y$ as arbitrary vector fields on $F$, we have on $F$
$$\left.d\left(\imath_{\hat A} \tom\right)(X,Y)\right|_p =\left.
X\tom(\hat
A,Y)\right|_p - \left.Y\tom(\hat A,X)\right|_p -
\left.\tom(\hat A,[X,Y])\right|_p = 0$$
by definition of $\tom$ along $F$. On the other hand
$$\left.d\left(\imath_{\hat A} \tom\right)(\hat B,Y)\right|_p =
\left.\hat
B\tom(\hat A,Y)\right|_p - \left.Y\tom(\hat
A,\hat B)\right|_p - \left.\tom(\hat A,[\hat B,Y])\right|_p=$$
$$= \left.\tom([\hat B,\hat A],Y)\right|_p - \left.Y\tom(\hat A,\hat
B)\right|_p =
\left.\tom([\hat B,\hat
A],Y)\right|_p + \left.d\mu_{\omega_\rho}(Y)(\tau([A,B]_\k))\right|_p = $$
$$= \left.\tom([\hat B,\hat A],Y)\right|_p +
\left.\omega_\rho(Y,\widehat{[A,B]_\k})\right|_p
= \left.\tom(\widehat{[A,B]_\k},Y)\right|_p -
\left.\tom(\widehat{[A,B]_\k},Y)\right|_p
= 0\ .$$
\par
\bigskip
Now,
for any sufficiently small $\epsilon\in \R^+$, we may
consider the
  $G$-invariant closed two-form $\omega_\epsilon$ on $M$
given
by
$$\omega_\epsilon = \pi^*\omega + \epsilon\
\tom\ ,\eqno(5.2)$$
By Lemma \ref{recovering} and (5.1), the restriction to
$F$ of the algebraic representative of $\omega_\epsilon$
is
$$Z_{\omega_\epsilon} = \epsilon \sum_{i = 1}^\co f_i Z_i
+ Z_\omega\ ,\qquad
\text{where}\qquad f_i \=
\mu_{\omega_\rho}(\tau(Z_i))\ .\eqno(5.3)$$
Using Proposition \ref{ricci}, we get that the restriction
to $F_{\operatorname{reg}}$ of the
algebraic representative of the Ricci form $\rho_\epsilon$
of
$\omega_\epsilon$ is given by
$$Z_{\rho_\epsilon} = \sum_{i=1}^\co \left(\phi_i +
\psi_{i\epsilon}
\right)Z_i + Z_V\ ,\eqno(5.4)$$ where
$$\phi_i \= \frac{J \hat Z_i(\det(f_{a,
b}))}{4 \pi \det(f_{a,b})} \ ,\qquad
\psi_{i\epsilon} \= \frac{\epsilon}{4 \pi}\sum_{\alpha\in
R^+_{\m}}\frac{a^j_{\alpha}\ f_{j,i}}{\epsilon\
a^j_{\alpha}\ f_{j}
+b_\alpha}\ ,\eqno(5.5)$$ and
$$f_{i,j} \= J \hat Z_j(f_i) = \omega_\rho(J \hat Z_j,
\hat Z_i)\ , \qquad a^j_{\alpha} \=\alpha(i Z_j)\ , \qquad
b_\alpha \=
\alpha(iZ_V)\ .\eqno(5.6)$$
We now notice that the map $Z_{\psi_\epsilon}|_F\=
\sum_{i=1}^\co \psi_{i \epsilon} Z_i$
defines a smooth closed $G$-invariant 2-form
$\psi_\epsilon$ on the regular part of $M$ by means
of (3.1) and (3.5); we claim that $\psi_\epsilon$ extends to a smooth
global 2-form on $M$, which is cohomologous to $0$. In
fact,
$Z_{\psi_\epsilon}|_F$ can be written as
$$Z_{\psi_\epsilon}|_F =\sum_i J \hat Z_i( \tilde
f_\epsilon) Z_i\ ,\qquad \text{with}\ \ \tilde f_\epsilon
\=
\frac{1}{4 \pi }\log\left(\prod_{\alpha \in R^+_\m}
\left(\epsilon\ a^j_\alpha\ f_{j} +
b_\alpha\right)\right)\ .$$
By the fact that $b_\alpha>0$ for every $\alpha\in
R_{\m}^+$, if $\epsilon$ is sufficiently small,
  the function $\tilde f_\epsilon$ is a well-defined
$K$-invariant function on $F$ and it extends to a
$G$-invariant
function on $M$, which we still denote by $\tilde
f_\epsilon$. Therefore, by
Lemma \ref{recovering}(b), it follows that $\psi_\epsilon$
coincides with the globally defined, two-form $- d d^c
\tilde f_\epsilon$
for any $\epsilon$ sufficiently small. \par
  From this we immediately get also
  that, for $\epsilon$ small, the Ricci form
$\rho_\epsilon$ is cohomologous to the two-form
$\rho_o \in c_1(M)$ given by
$$\rho_o = \rho_\epsilon - \psi_\epsilon = \rho_\epsilon +
d d^c \tilde f_\epsilon\ ,$$
whose algebraic representative on $F_{\operatorname{reg}}$ is equal
to
$$Z_{\rho_o} = \sum_{j=1}^\co\frac{J \hat Z_j(\det(f_{a,
b}))}{4 \pi \det(f_{a, b})} Z_j + Z_V\ .\eqno(5.7)$$
\par
\medskip
From $(5.6)_1$ and (4.4) , we may notice that
$$\left.\frac{J \hat Z_j(\det(f_{a, b}))}{4 \pi \det(f_{a,
b})}\right|_{F_{\reg}}= \frac{1}{4 \pi} \div(J\widehat{\tau(Z_i)})\
.\eqno(5.8)$$
Consider a basis $(W_1$, $\dots$, $W_\co)$ for $\t$ and a
 basis $(Z_1,\ldots,Z_\co,\ldots,Z_N)$ of
$\z(\k)$, which is $\B$-orthonormal and extends the set
$(Z_1,\ldots,Z_\co)$; Let
also
$(W_1^*$, $ \dots$, $W_\co^*)$ and
$(Z_1^*$, $\dots$, $Z_N^*)$ be the corresponding dual bases of
$\t^*$ and $\z(\k)^*$, respectively. We set $C = [c^i_j]$
to be the matrix defined by $\tau(Z_j) = \sum_i c^i_j
W_i$ with $c^i_j =0$ for $j\geq \co+1$ and
observe that $\tau^*(W_j^*) = \sum_{\ell=1}^\co c^j_\ell
Z_\ell^*$. Then, by (5.8),
we get that on $F$
$$- \left.\B\left( \sum_{j=1}^\co\frac{J \hat
Z_j(\det(f_{a, b})}{4 \pi \det(f_{a, b})} Z_j,
\cdot\right)\right|_{\z(\k)}=
\sum_{j=1}^\co\frac{J \hat Z_j(\det(f_{a,
b}))}{ 4 \pi \det(f_{a, b})}Z_j^* = $$
$$ = \sum_{j, \ell =1}^\co c^\ell_j \frac{J \hat
W_\ell(\det(f_{a, b}))}{4 \pi \det(f_{a, b})}Z_j^* =
  \frac{1}{4 \pi } \sum_{\ell=1}^\co
\div(J\widehat{W_\ell})\tau^*(W_\ell^*) = \tau^*\mu_\rho\
,\eqno(5.9)$$
  where the last equality is meaningful whenever $\rho$
is non-degenerate.
  Using Lemma \ref{recovering} (a) and (4.4), we see that
on $F_{\reg}$
$$\left.\rho_o(J\hat Z_i, \hat Z_j)\right|_{F_{\reg}} = J\hat
Z_i\left(\frac{J\hat
Z_j(\det(f_{a,b}))}{4 \pi \det(f_{a,b})}\right) =
\rho(J\widehat {\tau(Z_i)}, \widehat
{\tau(Z_j)})\eqno(5.10)$$
so that $\rho_o|_{TF} = \rho$.
Moreover, in case $\rho$ is non-degenerate, we also have
$$\ch Z_{\rho_o}|_F= \ch Z_V + \tau^* \mu_\rho\
.\eqno(5.11)$$
  \par
  Let us now prove the sufficiency of conditions in
Theorem 1.1. Assume that (1.1) holds and that $F$ is Fano
with $\rho >0$. We want to show that $\rho_o > 0$.
Indeed, from (5.10) and
the fact that the vector fields $\hat A$, $A\in \m$, are
$\rho_o$-orthogonal to $TF$
at all points of the fiber, we have that $\rho_o$
is positive if and only if the matrix $\left(\rho_o(\hat
F_\alpha,J\hat
F_\beta)\right)$ is positive definite at every point of $F
= \pi^{-1}(eK)$. We
now observe that the functions
$\rho_o(\hat F_{\alpha}, J\hat F_{\beta})$ vanish if
$\alpha\neq
\beta$, so
that we are reduced to check that
$$\rho_o(\hat F_{\alpha}, J\hat F_{\alpha}) =i
\alpha(Z_{\rho_o}) > 0\eqno(5.12)$$
for any $\alpha\in R_{\m}^+$ and at any point of the fiber
$F$. From (5.11), this turns out
to be
equivalent to (1.1). \par
Now, let us assume that $c_1(M)> 0$ and
let $\rho_1 \in c_1(M)$ be $G$-invariant and positive.
Being $\rho_1$ cohomologous to $\rho_o$ and
by (5.10), we have on $F$
$$\rho_1(\hat Z_i, J \hat Z_j)|_F = \rho_o(\hat Z_i, J\hat
Z_j)|_F +
d d^c\phi(\hat Z_i, J \hat Z_j)|_F = $$
$$ = \rho(\widehat{\tau(Z_i)}, J\widehat{ \tau(Z_j)}) +
d d^c\left(\phi|_F\right)(\widehat{\tau(Z_i)}, J\widehat{
\tau(Z_j)}) $$
for some $G$-invariant function $\phi$ on $M$. From this
it follows that $\rho_1|_{TF}$ is a positive 2-form in $c_1(F)$.
\par
Now, let us assume $\rho$ to be positive,
so that the last equality in (5.9) is meaningful. Recall
that, by Lemma \ref{recovering} (b), the algebraic
representatives of $\rho_1$ and $\rho_o$, restricted to
$F$, differ
  by a map $Z_\phi$ so that $Z_\phi|_F = - \sum_i J \hat
Z_i(\phi) Z_i$, for
some $K$-invariant smooth function $\phi: F \to \R$. In
particular, at the $T^\co$-fixed points of $F$,
the algebraic representatives of $\rho_1$ and $\rho_o$
coincides. On the other hand, we note that
$Z_{\rho_1}|_F$ is the $(-\B)$-dual of a moment map for the
action of $Z(K)$ on $F$
w.r.t. $\rho_1|_{TF}$, and therefore
$Z_{\rho_1}(F)\subset \z(\k)$ is a convex polytope with
vertices given by $Z_{\rho_1}(F^{Z(K)})$;
By the previous remark we see that
$Z_{\rho_1}(F)=Z_{\rho_o}(F)$.
By the fact that for any $\alpha \in R^+_\m$
$$ \rho_1(\hat X_\alpha, J \hat X_\alpha)|_F =
i\alpha(Z_{\rho_1}|_F)\ , $$
we have that $Z_{\rho_1}(F) =Z_{\rho_o}(F)\subset \CC$
and hence that $\ch Z_{\rho_o}(F) = \ch Z_V +\Delta_{\tau,
F} \subset \ch \CC$.
\bigskip
\bigskip
\section{Examples}
\bigskip
\noindent 1.\ The Hirzebruch surfaces $F_n$ ($n\in \mathbb N$) exhaust all
the homogeneous
toric bundles $M$ over flag manifolds when $\dim_\C M = 2$. The surface
$F_n$ can be
realized  as the homogeneous bundle
$\SL(2,\C)\times_{B,\tau}\CP^1$ over
$V = \CP^1 = \SL(2,\C)/B$. Here $B$ is the standard Borel subgroup of
$\SL(2,\C)$ given by
$B =\left\{\left(\begin{matrix} \alpha&\beta\cr
0&\alpha^{-1}\end{matrix}\right) \in \SL(2,\C)\right\}$
and $\tau:B\to \C^*$ is given by $\tau(\left(\begin{matrix}
\alpha&\beta\cr 0&\alpha^{-1}\end{matrix}\right)) = \alpha^n$, where
$\C^*$
acts on $\CP^1$ by $\zeta\mapsto \alpha\zeta$ for $\alpha\in\C^*$ and
$\zeta\in \C \cup \{\infty\}$.\par
It is well known that $F_n$ is Fano if and only if $n=0,1$ (see
\cite{Be1}). This property can be very rapidly established also by means of our
Theorem \ref{maintheorem}. In fact, using our notation and identifying
$\z(\k) = \k$ with $\R$ by means of the map  $\R\ni x\mapsto
\operatorname{diag}(ix,-ix)\in \su(2)$,  we have that
 $Z_V = -\frac{1}{8 \pi}$ and $\mathcal C =
\{x<0\}$. Then,
it is quite immediate to check that $-\B^{-1}(\Delta_{\tau,\CP^1}) =
[-\frac{n}{16 \pi },\frac{n}{16 \pi}]$,  so  that
$F_n$ is Fano if and only if $-\frac{1}{8 \pi } + [-\frac{n}{16 \pi },\frac{n}{16 \pi}]
\subset \{x<0\}$, i.e.  $n<2$.\par
\medskip
\noindent 2.\ Let us now assume that $F = \CP^\co$ and that $T^\co$ is the
standard maximal torus of $\SU(\co+1)$,
so that $(T^\co)^\C$ coincides with the
subgroup of diagonal matrixes in $\SL(\co+1,\C)$. Let us
also denote by $C= [c_i^j]$ the matrix with $\det C\neq
0$ with entries $c^i_j$ defined by $\tau(Z_i) =
\sum_{j=1}^\co c_i^j W_j$,
where the $W_j$ are the matrices in $\t$ defined by
$$W_j\= \frac{1}{\co+1}\cdot \operatorname{diag}(-i
,\ldots,\underset{\text{(j+1)-th place}}{\co i},\ldots,
-i) \in \t \subset \su(\co+1)\ .$$
In this case, the canonical polytope $\Delta_F$
  is the convex polytope with
vertices
$$Q_o = - \frac{1}{2\pi} \sum_{j = 1}^\co W_j^*\ , \qquad Q_r = Q_o +
\frac{\co + 1}{2 \pi} W^*_r\ ,\qquad 1 \leq r\leq \co\ .\eqno(6.1)$$
So, condition (1.1) amounts to say that all the
points $\ch P_o = \ch Z_V - \frac{1}{2 \pi}ÿ\sum_{i,j=1}^\co c_i^j \ch Z_i$ and
$\ch P_r = \ch P_o +
\frac{\co + 1}{2 \pi} \sum_{j =1 }^\co c_j^r \ch Z_j$, $1\leq r\leq
\co$, are in $\ch \CC$, or, equivalently, that the
points
$$P_o = Z_V - \frac{1}{2 \pi} \sum_{i,j=1}^\co c_i^j Z_i\ , \quad
P_r = P_o + \frac{\co+1}{2 \pi}ÿ \sum_{j=1}^\co c^r_jZ_j\ ,\quad
1 \leq r\leq \co\ ,\eqno(6.2)$$
are all in $\mathcal C$.\par
\medskip
Let us now construct explicitly an homogeneous toric
bundle with $c_1(M)>0$ and fiber
$F = \CP^2$.
  Consider the classical group $G=\SO(4n)$ with $\mathcal
B(X,Y) =
2(2n-1)\operatorname{Tr}(XY)$; we denote by $h_i$,
$1 \leq i \leq 2n$, the elements of
$\g=\so(4n)$ given by $h_i = E_{2i,2i+1} - E_{2i+1,2i}$,
where $E_{ij}$ denotes the matrix
whose unique non trivial entry is $1$ and in position $(i,j)$.
Given $J_1 = \sum_{i=1}^n h_i$ and
$J_2 = \sum_{i=n+1}^{2n}h_i$, we may consider the flag
manifold $G/K = \Ad(G)\cdot (J_1+2J_2)\cong
\SO(4n)/\U(n)\times \U(n)$. Using standard notations for
the roots, we have that
$$R_{\m}^+ = \{ \omega_i+\omega_j,\ 1\leq i <j \leq
2n\}\cup \{ \omega_j - \omega_i,\ 1\leq i \leq n < j \leq
2n\}\ .$$
We may consider $Z_i = \frac{J_i}{2\sqrt{n(2n-1)}}$ and
the homomorphism
$\tau:K\to \SU(3)$ with $\tau|_{K_{ss}} = e$ and
$\tau(Z_i) = c_i^ jW_j$ with $c_i^ j = \frac{3
n\delta_i^j}{2\sqrt{n(2n-1)}}$.
Since $Z_2\in {\mathcal C}$, by (6.2) we have that the manifold
$M\=G\times_{K,\tau}\CP^2$
is Fano if and only if
$$P_o = Z_V - \frac{1}{2 \pi}ÿ\sum_{i,j} c_i^ jZ_i = Z_V -
 \frac{1}{2 \pi}  \frac{3 n}{\sqrt{4n(2n-1)}}(Z_1 + Z_2) = $$
$$ =\frac{1}{8 \pi  (2n -1)}\left(\sum_{1\leq i < j \leq 2n} ( h_i + h_j) +
\sum_{\smallmatrix1 \leq i \leq n
\\ n +1 \leq i \leq 2n
\endsmallmatrix} (h_j - h_i) - 3 \sum_{i = 1}^{2n} h_i\right)$$
 and
 $$P_1 = P_o +\frac{9 n}{2 \pi \sqrt{4 n(2n-1)}}Z_1 = P_o +
\frac{9}{8 \pi (2n-1)}\sum_{i = 1}^{n} h_i $$
 are both in $\mathcal C$.
A direct inspection shows that this occurs if and only if
$n \geq 5$.

\bigskip
\bigskip
\font\smallsmc = cmcsc8
\font\smalltt = cmtt8
\font\smallit = cmti8
\hbox{\parindent=0pt\parskip=0pt
\vbox{\baselineskip 9.5 pt \hsize=3.1truein
\obeylines
{\smallsmc
Fabio Podest\`a
Dip. Matematica e Appl. per l'Architettura
Universit\`a di Firenze
Piazza Ghiberti 27
I-50100 Firenze
ITALY
}\medskip
{\smallit E-mail}\/: {\smalltt podesta@math.unifi.it
}
}
\hskip 0.0truecm
\vbox{\baselineskip 9.5 pt \hsize=3.7truein
\obeylines
{\smallsmc
Andrea Spiro
Dip. Matematica e Informatica
Universit\`a di Camerino
Via Madonna delle Carceri
I-62032 Camerino (Macerata)
ITALY
}\medskip
{\smallit E-mail}\/: {\smalltt andrea.spiro@unicam.it}
}
}

\end{document}